\documentclass[10pt,twocolumn]{article}

\usepackage[utf8]{inputenc}
\usepackage{amsmath}
\usepackage{amssymb}
\usepackage{amsthm}
\usepackage{titlesec}
\usepackage{cite}
\usepackage{color}
\usepackage{booktabs}
\usepackage{graphicx}
\usepackage{nicefrac}
\usepackage[colorlinks=false]{hyperref}
\usepackage[format=plain,labelfont=it]{caption}
\usepackage[left=1.5cm,right=1.5cm,top=2cm,bottom=2cm]{geometry}

\titleformat{\section}{\centering\normalfont\scshape}{\Roman{section}.}{5pt}{}
\titleformat{\subsection}{\normalfont\it}{\Alph{subsection}.}{5pt}{}
\titleformat{\subsubsection}{\normalfont\it}{\hspace{4mm}\arabic{subsubsection})}{5pt}{}

\newcommand\infoFootnote[1]{%
  \begingroup
  \renewcommand\thefootnote{}\footnote{#1}%
  \addtocounter{footnote}{-1}%
  \endgroup}

\newtheorem{thm}{Theorem}
\newtheorem{cor}[thm]{Corollary}
\newtheorem{lem}[thm]{Lemma}

\newcommand{\R}{\mathbb{R}} 
\newcommand{\N}{\mathbb{N}} 

\newcommand{\Ac}{\mathcal{A}}
\newcommand{\Fc}{\mathcal{F}}
\newcommand{\Rc}{\mathcal{R}}
\newcommand{\Sc}{\mathcal{S}} 
\newcommand{\Tc}{\mathcal{T}} 
\newcommand{\Uc}{\mathcal{U}} 
\newcommand{\Xc}{\mathcal{X}} 

\newcommand{\Ab}{\boldsymbol{A}}
\newcommand{\Bb}{\boldsymbol{B}}

\newcommand{\Eb}{\boldsymbol{E}}
\newcommand{\Fb}{\boldsymbol{F}}
\newcommand{\Gb}{\boldsymbol{G}}
\newcommand{\Hb}{\boldsymbol{H}}
\newcommand{\Ib}{\boldsymbol{I}}

\newcommand{\Kb}{\boldsymbol{K}}

\newcommand{\Pb}{\boldsymbol{P}}
\newcommand{\Qb}{\boldsymbol{Q}}
\newcommand{\Rb}{\boldsymbol{R}}
\newcommand{\Sb}{\boldsymbol{S}}

\newcommand{\Ub}{\boldsymbol{U}}

\newcommand{\bb}{\boldsymbol{b}}
\newcommand{\cb}{\boldsymbol{c}}
\newcommand{\db}{\boldsymbol{d}}
\newcommand{\eb}{\boldsymbol{e}}
\newcommand{\fb}{\boldsymbol{f}}

\newcommand{\rb}{\boldsymbol{r}}

\newcommand{\ub}{\boldsymbol{u}}
\newcommand{\vb}{\boldsymbol{v}}

\newcommand{\xb}{\boldsymbol{x}}

\newcommand{\Omegab}{\boldsymbol{\Omega}}
\newcommand{\Thetab}{\boldsymbol{\Theta}}

\newcommand{\lambdab}{\boldsymbol{\lambda}}
\newcommand{\deltab}{\boldsymbol{\delta}}

\newcommand{\zerob}{\boldsymbol{0}}
\newcommand{\oneb}{\boldsymbol{1}}

\newcommand{\diag}{\mathrm{diag}}
\newcommand{\norm}[1]{\left\lVert#1\right\rVert}

\DeclareMathOperator*{\find}{find}

\title{\vspace{-2mm}\bf Efficient computation of Lipschitz constants for MPC with symmetries}
\author{Dieter Teichrib and Moritz Schulze Darup\vspace{2mm}}
\date{}

\begin{document}

\maketitle

\textbf{\textit{Abstract}.} {\bf Lipschitz constants for linear MPC are useful for certifying inherent robustness against unmodeled disturbances or robustness for neural network-based approximations of the control law. In both cases, knowing the minimum Lipschitz constant leads to less conservative certifications. Computing this minimum Lipschitz constant is trivial given the explicit MPC. However, the computation of the explicit MPC may be intractable for complex systems. The paper discusses a method for efficiently computing the minimum Lipschitz constant without using the explicit control law. The proposed method simplifies a recently presented mixed-integer linear program (MILP) that computes the minimum Lipschitz constant. The simplification is obtained by exploiting saturation and symmetries of the control law and irrelevant constraints of the optimal control problem.}
\infoFootnote{D. Teichrib and M. Schulze Darup are with the \href{https://rcs.mb.tu-dortmund.de/}{Control and~Cyber-physical Systems Group}, Faculty of Mechanical Engineering, TU Dortmund University, Germany. E-mails:  \href{mailto:moritz.schulzedarup@tu-dortmund.de}{\{dieter.teichrib, moritz.schulzedarup\}@tu-dortmund.de}. \vspace{0.5mm}}
\infoFootnote{\hspace{-1.5mm}$^\ast$This paper is a \textbf{preprint} of a contribution to the 62nd IEEE Conference on Decision and Control 2023.}

\section{Introduction and problem statement}

Model predictive control (MPC) is an established method for the performance-oriented control of dynamical systems subject to state and input constraints (see \cite{Rawlings2017} for an overview). For various applications or realizations of MPC, knowing a Lipschitz constant of the control law is beneficial. For instance, Lipschitz constants can be used to certify inherent robustness of classical MPC against unmodeled disturbances \cite{Scokaert1997,Limon2002}. Based on similar concepts, Lipschitz constants can be used to certify robustness for neural network-based approximations of MPC laws \cite{Fabiani2022,Teichrib2023}. 

As pointed out in \cite{SchulzeDarup2017}, computing a Lipschitz constant for linear MPC is trivial, if the (piecewise affine) control law is explicitly known. However, the number of affine segments may grow exponentially with the state dimension and the number of constraints in the optimal control problem (OCP). Thus, explicitly computing the control law often becomes intractable for complex systems. This motivates the design of methods that allow the computation of a Lipschitz constant without explicitly computing the control law. Such methods have been previously addressed, e.g, in \cite{SchulzeDarup2017,Fabiani2022}. In \cite{SchulzeDarup2017}, a procedure enumerating potential active sets of the OCP is proposed. However, the procedure is inefficient and, in parts, based on an unproven conjecture (see \cite[Conj.~4]{SchulzeDarup2017}).
In \cite{Fabiani2022}, the problem of finding the smallest Lipschitz constant for MPC is formulated as a mixed-integer linear program (MILP). The procedure is elegant and it can be applied in case where the computation of the explicit control law is numerically intractable (cf. \cite[Tab.~I]{Fabiani2022}). Still, solving the resulting MILP may be time-consuming since, typically, many binary variables are involved. 

In this paper, we aim for a more efficient computation of Lipschitz constants via MILP. To this end, we first slightly improve the MILP formulation from \cite{Fabiani2022} by reducing the initial number of binary variables. Afterwards, we present various preprocessing steps to further reduce the number of binary variables in the MILP. In this context, the most powerful reduction step builds on symmetries that often arise in MPC (see, e.g., \cite{Danielson2015}). Furthermore, the observation that the control law is constant in many regions for a large prediction horizon also allows a significant reduction in the number of binary variables.   

The paper is organized as follows. The notation is given in the remainder of this section. Section~\ref{sec:Preliminaries} is devoted to the basics of MPC and presents a known representation of the OCP in terms of an MILP, which can be used to compute the minimum Lipschitz constant of the control law. In Section~\ref{sec:ReducOfBinVar}, we describe several methods for reducing the number of binary variables needed to compute the minimum Lipschitz constant via MILP. We illustrate the effectiveness of the proposed method with some examples from the literature in Section~\ref{sec:Examples} and give a conclusion in Section~\ref{sec:Conclusion}.    

\vspace{2mm}
\textit{Notation.} For matrices $\Kb\in\R^{m \times n}$, we denote its $i$-th row by $\Kb_{i,:}$, its $j$-th column by $\Kb_{:,j}$, and its elements by $\Kb_{i,j}$, respectively. We further use the shorthand notation $\Kb_{i}$ for $\Kb_{i,:}$ and we indicate different instances of a matrix $\Kb$ by a superscript, e.g., $\Kb^{(i)}$. The product 
\begin{equation*}
\Kb \Xc := \{\Kb\xb\in\R^m \ | \ \xb \in \Xc \}
\end{equation*}
is defined for a compact and convex set $\Xc$. For a vector $\vb\in\R^m$ we define $\diag(\vb)$ as a diagonal matrix with the elements of $\vb$ being the diagonal elements. Moreover, $\Ib$ refers to the identity matrix and we denote column vectors or matrices whose entries are all $0$ respectively $1$ by $\zerob$ respectively~$\oneb$. Finally, all inequalities involving matrices or vectors are understood element-wise.

\section{Preliminaries}
\label{sec:Preliminaries}

\subsection{Lipschitz constants of piecewise affine functions}\label{SubSec:LpForPWA}
In general a Lipschitz constant of a function $\fb:\R^n \rightarrow \R^m$ on a domain $\Fc \subseteq \R^n$ is a  constant $L$ satisfying
\begin{equation*}
    \| \fb(\xb)-\fb(\tilde{\xb})\| \leq L \| \xb - \tilde{\xb}\|
\end{equation*}
for all $\xb,\tilde{\xb}\in\Fc$ and some vector norm $\|\cdot\|$. 
We here focus on $p$-norms and, in particular, the cases $p\in \{1,\infty\}$. We denote corresponding Lipschitz constants by $L_p$. Further, we aim for as small as possible Lipschitz constants and denote the smallest one by $L_p^\ast$. Now, under the assumption that $\fb$ is a continuous piecewise affine (PWA) function on $\Fc$, i.e., of the form
\begin{equation}
\label{eq:fPWA}
    \fb(\xb)= \left\{ \begin{array}{ll}
\Kb^{(1)} \xb + \bb^{(1)} & \text{if}\,\,\,\xb\in\Rc^{(1)},\\
\vdots \\
\Kb^{(s)} \xb + \bb^{(s)} & \text{if}\,\,\,\xb\in\Rc^{(s)},
\end{array}\right.
\end{equation}
with $\Kb^{(i)} \!\in \R^{m \times n}$ and $\bb^{(i)} \!\in \R^m$, then it is well-known~that 
\begin{equation}
\label{eq:maxKnorm}
L_p^\ast= \max_{i \in \{1,\dots, s\}} \| \Kb^{(i)} \|_p,
\end{equation}
where $\|\cdot \|_p$ here refers to the matrix norm induced by the vector $p$-norm (see, e.g., \cite[Prop.~3.4]{Gorokhovik1994}). 
As pointed out in \cite[Lem.~5.1]{Fabiani2022}, the matrix norms $\| \Kb \|_1$ and $\| \Kb \|_\infty$ can be evaluated by solving a linear program (LP).
In fact, one finds
\begin{align}
    \label{eq:normLP}
    \norm{\Kb}_1 &= \min_{l,\vb^{(1)},\dots,\vb^{(n)}}  l \\
    \nonumber
    \text{s.t.} \,\,\,\, &\oneb^\top \vb^{(j)} \leq l, \,\,\,\,-\vb^{(j)} \leq \Kb_{:,j} \leq \vb^{(j)}, \,\,\,\, \forall j \in \{1,\dots,n\}.
\end{align}
Due to $||\Kb||_{\infty}=||\Kb^\top||_{1}$, the LP for computing the $\infty$-norm can be formulated in a similar way.

\subsection{Linear MPC via MILP}
MPC for linear discrete-time systems builds on solving an OCP of the form
\begin{align}
    \label{eq:OCP}
    V(\xb) := \!\!\!\!\min_{\substack{\xb(0),...,\xb(N)\\ \ub(0),...,\ub(N-1)}} 
    \!\!\!\!\!\!\!\!\!\!\!\!\!\!\!\!\!\!\!\! & 
    \,\,\,\,\,\,\,\,\,\,\,\,\,\,\,\,\,\,\,\,\,
    \varphi( \xb(N)) + \! \sum_{k=0}^{N-1} \ell(\xb(k),\ub(k))   \\
    \nonumber
    \text{s.t.} \quad \quad  \xb(0)&=\xb, \\
    \nonumber
    \xb(k+1)&=\Ab\,\xb(k) + \Bb \ub(k), \quad\!\!\forall k \in \{0,...,N-1\}, \\
    \nonumber
    \left(\xb(k),\ub(k)\right) & \in \Xc \times \Uc, \quad\hspace{14.9mm}\forall k \in \{0,...,N-1\}, \\
    \nonumber
    \xb(N) & \in \Tc
\end{align}
in every time step for the current state $\xb$. 
Here, $N\in\N$ refers to the prediction horizon and 
\begin{equation}
\nonumber
    \varphi(\xb):= \xb^\top \Pb \xb \quad \text{and} \quad \ell(\xb,\ub):=\xb^\top \Qb \xb + \ub^\top \Rb \ub
\end{equation}
denote the terminal and stage cost, respectively, where we assume that the weighting matrices  $\Pb$, $\Qb$, and $\Rb$ are positive definite. The dynamics of the linear prediction model are described by $\Ab\in \R^{n \times n}$ and ${\Bb \in \R^{n\times m}}$. State and input constraints can be incorporated via the polyhedral sets $\Xc$ and $\Uc$. Finally, a polyhedral terminal set $\Tc$ allows to enforce closed-loop stability 
(see \cite{Mayne2000} for details).
Condensing the OCP~\eqref{eq:OCP} leads to a parametric quadratic program (QP) of the form 
\begin{align}
    \label{eq:QP}
    \Ub^\ast(\xb) := \arg &\min_{\Ub}  \frac{1}{2} \Ub^\top \Hb \Ub + \xb^\top \Fb^\top \Ub \\ 
   \nonumber
    & \text{s.t.} \quad \Gb \Ub \leq \Eb \xb + \db
\end{align}
with the decision variable $\Ub\in \R^{mN}$ reflecting the predicted input sequence (i.e, the stacked vectors $\ub(0),...,\ub(N-1)$) and with $\Hb,\Fb,\Gb,\Eb,\db$ denoting condensed matrices obtained from $\Ab,\Bb,\Pb,\Qb,\Rb$ and the specifications of $\Xc,\Uc,\Tc$ (see \cite[Chap.~3]{Maciejowski2002} for details). Now, MPC typically builds on applying the first element of the optimal input sequence, i.e., $\ub^\ast(0)$, and repeating the procedure at the next sampling instant. Hence, the resulting control law $\fb_{\text{MPC}}:\Fc_{\text{MPC}} \rightarrow \Uc$ can be defined as
\begin{equation}
    \label{eq:gMPC}
    \fb_{\text{MPC}}(\xb):=\Sb \,\Ub^\ast(\xb),
\end{equation}
where $\Sb:=(\Ib \ \zerob \ \dots \ \zerob) \in \R^{m \times mN}$ serves as a selection matrix and where the set $\Fc_{\text{MPC}}$ contains all $\xb \in \R^n$ for which \eqref{eq:QP} (or, equivalently, \eqref{eq:OCP}) is feasible. Remarkably, $\fb_{\text{MPC}}$ is of the form \eqref{eq:fPWA} with regions $\Rc^{(i)}$ representing polyhedral sets with pairwise disjoint interiors \cite[Thm.~4]{Bemporad2002}. 

Clearly, the solution of~\eqref{eq:QP} can be characterized by the Karush-Kuhn-Tucker (KKT) conditions 
\begin{subequations}
    \label{eq:KKT}
 \begin{align}
 \label{eq:KKT_opt}
    \Hb \Ub^*(\xb) + \Fb \xb + \Gb^\top \lambdab^*(\xb) &= \zerob, \\
    \label{eq:KKT_r}
   \rb^\ast(\xb) &= \Eb \xb + \db - \Gb \Ub^\ast(\xb), \\
   \label{eq:KKT_rGeq0}
    \rb^\ast(\xb) &\geq \zerob, \\
    \label{eq:KKT_lambdaGeq0}
    \lambdab^*(\xb) &\geq \zerob, \\
    \label{eq:KKT_complement}
    \diag(\lambdab^*(\xb) ) \,\rb^\ast(\xb) & = \zerob
\end{align}   
\end{subequations}
where $q\in \N$ reflects the dimension of $\db$, i.e., the number of constraints.
Now, let us assume upper bounds $\overline{\lambdab}$ and $\overline{\rb}$ for the dual optimizers $\lambdab^*(\xb)$ respectively the residuals $\rb^\ast(\xb)$ are known, i.e., $\lambdab^\ast(\xb)\leq \overline{\lambdab}$ and $\rb^\ast(\xb)\leq \overline{\rb}$ for all $\xb \in \Fc_{\text{MPC}}$. Since \eqref{eq:QP} is feasible for all $x \in \Fc_{\text{MPC}}$ by definition, such upper bounds exist and are finite. Then, the conditions~\eqref{eq:KKT_rGeq0}--\eqref{eq:KKT_complement} can be rewritten as
\begin{subequations}
\label{eq:KKT2MI}
 \begin{align}
 \label{eq:KKT2MI_r}
0 \leq \rb^\ast(\xb) &\leq \diag(\overline{\rb}) (\oneb-\deltab^\ast(\xb)), \\
0 \leq \lambdab^*(\xb) &\leq \diag(\overline{\lambdab}) \,\deltab^\ast(\xb), \\
 \deltab^\ast(\xb) &\in \{0,1\}^q
\end{align} 
\end{subequations}
\cite[Eq.~(13)]{Fabiani2022}. Thus, the QP \eqref{eq:QP} can be solved by solving the MI feasibility problem 
\begin{align}
    \label{eq:QP_MIFeas}
    \text{find}\,\,\,\Ub^*,\,\lambdab^*,\,\rb^\ast,\,\text{and}\,\,\deltab^\ast \quad \text{s.t.} \quad  \text{\eqref{eq:KKT_opt}--\eqref{eq:KKT_r} and \eqref{eq:KKT2MI}},
\end{align}
where we omit dependencies on $\xb$ for brevity.

\subsection{Local MPC gain via MILP}\label{sec:KMPCviaMILP}

Given a solution to~\eqref{eq:QP_MIFeas}, it is easy to see that the set
$$
\Ac(\xb):=\{ i \in \{1,\dots,q\}\,|\, \deltab_i^\ast(\xb)=1\}
$$
reflects the active constraints for~\eqref{eq:QP}. In principle, $\Ac(\xb)$ allows to compute the local MPC segment, i.e., $\Kb^{\ast}$, $\bb^\ast$, and $\Rc^\ast$ such that $\xb \in \Rc^\ast$ and 
$$
\fb_{\text{MPC}}(\tilde{\xb}) = \Kb^\ast \tilde{\xb} + \bb^\ast \quad \forall\,\, \tilde{\xb} \in \Rc^\ast.
$$
However, performing the computation analogously to \cite{Bemporad2002} leads to nonlinear relations between $\Kb^{\ast}$ and $\deltab^\ast(\xb)$. Hence, aiming for a combination of \eqref{eq:maxKnorm}, \eqref{eq:normLP} and \eqref{eq:QP_MIFeas} in one MILP, we need to derive $\Kb^\ast$ differently. A suitable approach has been proposed in the proof of \cite[Thm.~5.3]{Fabiani2022}. 
The underlying idea is to sample the MPC segment (potentially even outside its domain~$\Rc^\ast$) and to use the samples to characterize $\Kb^\ast$. The sampling points are chosen as
\begin{equation*}
\xb^{(j)} := \xb + \eb^{(j)} \quad \forall j\in\{1,\dots,n\},
\end{equation*}
where $\eb^{(j)}$ refers to the $j$-th canonical unit vector. One then constructs corresponding $\Ub^{(j)}$ from the MPC segment around $\xb$ by introducing the additional variables $\Ub^{(j)}$, $\lambdab^{(j)}$, and $\rb^{(j)}$ and by  augmenting the constraints \eqref{eq:KKT_opt}--\eqref{eq:KKT_r} and \eqref{eq:KKT2MI} with 
\begin{subequations}
\label{eq:augmentedConstraints}
\begin{align}
        \!\Hb \Ub^{(j)}\! + \Fb (\xb + \eb^{(j)}) + \Gb^\top \lambdab^{(j)} &= \zerob \\
        \Eb (\xb + \eb^{(j)}) + \db - \Gb \Ub^{(j)} &= \rb^{(j)}  \\
        -M (\oneb-\deltab^*(\xb)) \leq \rb^{(j)} &\leq M (\oneb-\deltab^*(\xb))  \\
        \label{eq:KKT_MI_Lambda_j}
        -M \deltab^*(\xb) \leq \lambdab^{(j)} &\leq M \deltab^*(\xb) 
    \end{align}
\end{subequations}
for all $j \in \{1,\dots,n\}$, where $M$ is a sufficiently large number such that feasible solutions satisfy ${\rb_i^{(j)},\lambdab_i^{(j)} \!\in (-M,M)}$.
Then, we obviously have $\Sb\big(\Ub^{(j)}-\Ub^\ast(\xb)\big) = \Kb^\ast \eb^{(j)}$ and consequently
\begin{equation}
    \label{eq:KfromUs}
\Kb^\ast = \Sb \begin{pmatrix}
\Ub^{(1)}-\Ub^\ast(\xb) & \dots & \Ub^{(n)}-\Ub^\ast(\xb)  
\end{pmatrix}.
\end{equation}
We could now compute $\|\Kb^\ast\|_1$ by solving \eqref{eq:normLP} for $\Kb^\ast$ as in \eqref{eq:KfromUs}  subject to the additional constraints \eqref{eq:KKT_opt}--\eqref{eq:KKT_r}, \eqref{eq:KKT2MI}, and \eqref{eq:augmentedConstraints}, which results in an MILP. However, including a maximization analogously to~\eqref{eq:maxKnorm} (here over $\xb \in \Xc$) is non-trivial since \eqref{eq:normLP} calls for a minimization. To circumvent this issue, one can consider the dual of \eqref{eq:normLP} and observe that the corresponding optimizer is binary \cite[Lem.~5.1]{Fabiani2022}. This finally allows to compute the Lipschitz constant $L_1^\ast$ for an MPC scheme based on an MILP with 
\begin{equation}
\label{eq:numBinaryFabiani1}
  q+(2m+1)n  
\end{equation}
binary variables (or $q+(2n+1)m$ for $L_\infty^\ast$) \cite{Fabiani2022}. 

\section{Reducing binary variables for efficiency}\label{sec:ReducOfBinVar}

The numerical complexity for solving an MILP crucially depends on the number of binary variables. Hence, in order to enable more efficient Lipschitz constant computations for MPC, we aim for a reduction of the number of binary variables in the corresponding MILP. In this context, we first stress that the number in~\eqref{eq:numBinaryFabiani1} consists of two terms resulting from the $q$ constraints of the MPC-related QP in~\eqref{eq:QP} and the $(2m+1)n$ constraints in~\eqref{eq:normLP}. Hence, we have two immediate options for reducing the number of binary variables: First, reducing the number of constraints in~\eqref{eq:QP} relevant for the computation of  Lipschitz constants and, second, implementing the norm evaluation more efficiently. 
We investigate these options in Sections~\ref{subsec:moreEfficientNorm} and \ref{sec:ExcMPCConstraints}, respectively.
Finally, we show in Section \ref{sec:symmetries} that exploiting symmetries, which are often present in MPC, can be beneficial for the computation of Lipschitz constants.

\subsection{More efficient norm computation}
\label{subsec:moreEfficientNorm}

Roughly speaking, the computation of Lipschitz constants  proposed in \cite{Fabiani2022} builds on reformulations of the QP~\eqref{eq:QP} and the LP~\eqref{eq:normLP} in terms of MI feasibility problems.
In \cite{Fabiani2022}, the reformulation of~\eqref{eq:normLP} is realized based on the LP's dual. 
Next, we propose a direct reformulation based on the primal LP, 
which is inspired by MI modeling techniques from, e.g., \cite{Fischetti2018} and which requires fewer binary variables than~\cite{Fabiani2022}. 

\begin{lem}
\label{lem:normMIFeas}
    Let $\Kb \in \R^{m\times n}$ and consider the  conditions
\begin{subequations}
\label{eq:newNormConditions}
\begin{align}
\label{eq:newNormConditionsPlusMinus}
        \Kb &=\Kb^+ - \Kb^-, \\
        \label{eq:newNormConditionsPlusPos}
        \zerob \leq \Kb^{+} &\leq M \big(
           \, \deltab^{(1)} \,\,\,\, \dots \,\,\,\,   \deltab^{(m)} \,\big)^\top,\\ 
        \label{eq:newNormConditionsMinusPos}
        \zerob \leq  \Kb^{-} &\leq M \Big(\oneb- \big(
            \deltab^{(1)} \,\,\,\, \dots \,\,\,\,   \deltab^{(m)} \big)^\top \Big), \\
        \label{eq:newNormConditionsC}
        \cb^\top &=\oneb^\top (\Kb^+ + \Kb^-), \span \span\\
        \label{eq:newNormConditionsAlpha}
        \cb \leq  \oneb l &\leq \cb + M(\oneb-\deltab^{(m+1)}), \\
        \label{eq:newNormConditionsOnlyOne}
        \oneb^\top \deltab^{(m+1)} &= 1, \\
        \deltab^{(1)},\dots,\deltab^{(m+1)} &\in \{0,1\}^n
    \end{align}
    \end{subequations}
    for some $M$ being larger than the largest absolute value of the entries in $\Kb$. Then, any solution to 
   \begin{equation}
   \label{eq:newMIFeasNorm}
        \find \ l,\,  \Kb^+,\, \Kb^-,\,\cb,\,\deltab^{(1)},\dots,\deltab^{(m+1)} \quad \text{s.t.} \quad \eqref{eq:newNormConditions} \quad\,
   \end{equation}
    is such that $\|\Kb\|_1=l$.
\end{lem}

\begin{proof}
    It is easy to see that the conditions \eqref{eq:newNormConditionsPlusMinus}--\eqref{eq:newNormConditionsMinusPos} together with $\deltab^{(1)},\dots,\deltab^{(m)}\in \{0,1\}^n$ imply
$$
(\Kb_{i,j}^+,\Kb_{i,j}^-) := \left\{ \begin{array}{ll}
(\Kb_{i,j},0) & \text{if} \quad \Kb_{i,j}\geq 0, \\
(0,-\Kb_{i,j}) & \text{otherwise}.
\end{array}\right. 
$$
As a consequence, the entries of $\cb$ as in~\eqref{eq:newNormConditionsC} reflect the absolute column sums of $\Kb$.  Due to \eqref{eq:newNormConditionsAlpha}, \eqref{eq:newNormConditionsOnlyOne}, and $\deltab^{(m+1)} \in \{0,1\}^n$, $l$ equals the largest entry of $\cb$, which is, by definition, identical to $\|\Kb\|_1$.
\end{proof}
\vspace{2mm}

Obviously, the MI feasibility problem~\eqref{eq:newMIFeasNorm} only requires $(m+1)n$ binary variables and thus $mn$ less than the counterpart in \cite{Fabiani2022}. Again, a similar problem can easily be constructed to compute $\|\Kb\|_\infty$ using $(n+1)m$ binary variables. Likewise, the reduction compared to \cite{Fabiani2022} amounts~to~$mn$.

\subsection{Excluding MPC constraints}
\label{sec:ExcMPCConstraints}

We propose two approaches to reduce the number of constraints in~\eqref{eq:QP} and consequently the number of binary variables in \eqref{eq:QP_MIFeas}.
The approaches are conceptually decoupled but it will turn out that their implementation can be efficiently coupled.
The first approach builds on the straightforward observation that the $i$-th constraint in~\eqref{eq:QP} (i.e., ${\Gb_i \Ub \leq \Eb_i \xb + \db_i}$) is irrelevant for the MPC scheme, if there exists no $\xb \in \Fc_{\text{MPC}}$ such that $\Gb_i\Ub^*(\xb)=\Eb_i\xb + \db_i$.  
According to the following lemma, such a situation can be identified based on an MI feasibility problem similar to \eqref{eq:QP_MIFeas}.

\begin{lem}
\label{lem:iConstraintInfeasible}
    Let $i \in \{1,\dots,q\}$. If the MI feasibility problem \eqref{eq:QP_MIFeas} with the additional decision variable $\xb$ and with the additional constraint $\deltab_i^\ast=1$ is infeasible, then 
    \begin{equation*}
    \Gb_i\Ub^*(\xb)-\Eb_i\xb <\db_i \quad  \text{for all} \ \xb \in \Fc_{\text{MPC}}.
    \end{equation*}
\end{lem}
\vspace{2mm}
\begin{proof}
We first note that the unmodified problem \eqref{eq:QP_MIFeas} has, by construction, the same feasible set $\Fc_{\text{MPC}}$ as~\eqref{eq:QP} with respect to the parameter $\xb$. Now, the additional constraint $\deltab_i^\ast=1$ implies $\rb_i^\ast=0$ in \eqref{eq:KKT2MI_r} and consequently $\Gb_i \Ub^\ast = \Eb_i \xb + \db_i$ in \eqref{eq:KKT_r}. Hence, 
infeasibility of the MI feasibility problem in the claim with $\xb$ as a decision variable immmediately implies that there exists no $\xb \in \Fc_{\text{MPC}}$ such that $\Gb_i\Ub^*(\xb)-\Eb_i\xb = \db_i$. In other words, solving \eqref{eq:QP} for any feasible $\xb \in \Fc_{\text{MPC}}$ results in $\Gb_i\Ub^*(\xb)-\Eb_i\xb <\db_i$.
\end{proof}
\vspace{2mm}

Clearly, with the help of Lemma~\ref{lem:iConstraintInfeasible}, we can eliminate irrelevant constraints by checking the corresponding MI feasibility problem for every (or some) $i\in\{1,\dots,q\}$. At this point, it might seem counterintuitive to consider multiple MI problems in order to simply the overlaying MILP of interest. However, our numerical benchmark in Section~\ref{sec:Examples} clearly shows that this approach is meaningfull and that the overall runtime can be (significantly) shortened compared to a direct solution of the unmodified MILP for the computation of Lipschitz constants. Remarkably, Lemma~\ref{lem:iConstraintInfeasible} could also be used to simplify an MPC scheme offline in order to accelerate the QP solutions online.

The second approach for reducing binary variables associated to constraints differs from the first one in that it is tailored to the problem at hand. It is based on the observation that the MPC law of the form \eqref{eq:fPWA} often contains many segments with $\Kb^{(i)}=\zerob$ resulting in constant inputs determined by the bias term $\bb^{(i)}$. In particular, this situation often arises if box-shaped input constraints are present. Clearly, due to $\|\zerob\|_p=0$, the constant segments are irrelevant for the computation of Lipschitz constants (but they obviously matter for the MPC scheme). Now, according to the following theorem, our procedure to identify and exclude some of these segments is similar to the first approach.

\begin{thm}
\label{thm:saturatingInputs}
Let $i \in \{1,\dots,q\}$ and let $\Fc_{\text{MPC}}$ be full-dimensional. Consider the MILP
\begin{equation}
\label{eq:costDeltaU}
 \max_{\substack{\Ub^{\ast}, \Ub^{(1)},\dots,\Ub^{(m)}, \xb, \xb^{(1)},\dots,\xb^{(m)},\\ \lambdab^\ast, \lambdab^{(1)}, \dots,\lambdab^{(m)}, \rb^\ast, \rb^{(1)},\dots, \rb^{(m)},\deltab^\ast}} \sum_{j=1}^m \Ub_j^{(j)}-\Ub_j^\ast 
\end{equation}
subject to \eqref{eq:KKT_opt}--\eqref{eq:KKT_r}, \eqref{eq:KKT2MI}, $\deltab_i^\ast=1$, and
\begin{subequations}
\label{eq:addConstraintsSaturation}
   \begin{align}
        \Hb \Ub^{(j)}+ \Fb \xb^{(j)} + \Gb^\top \lambdab^{(j)} &= \zerob, \\
        \Eb \xb^{(j)} + \db - \Gb \Ub^{(j)} &= \rb^{(j)},  \\
        \zerob \leq \rb^{(j)} &\leq \diag(\overline{\rb}) (\oneb-\deltab^\ast),  \\
\zerob \leq \lambdab^{(j)} &\leq \diag(\overline{\lambdab}) \,\deltab^\ast
    \end{align} 
\end{subequations}
for every $j\in\{1,\dots,m\}$. If the MILP is feasible and
returns $0$ as the optimal objective function value, then the $i$-th constraint can be omitted for the computation of the Lipschitz constant.
\end{thm}

\begin{proof}
We initially neglect the objective function in~\eqref{eq:costDeltaU} and investigate the corresponding MI feasibility problem. 
We further assume feasibility since the theorem is irrelevant otherwise. Now, we consider any feasible set of decision variables and note that, based on the corresponding $\Ub^\ast$, $\xb$, $\lambdab^\ast$, $\rb^\ast$, and $\deltab^\ast$, we can construct another set of feasible variables by choosing $\Ub^{(j)}:=\Ub^\ast$,  $\xb^{(j)}:=\xb$, $\lambdab^{(i)}:=\lambdab^\ast$, and $\rb^{(j)}:=\rb^\ast$ for every $j\in\{1,\dots,m\}$ (and keeping the other variables). Clearly, the associated objective function value is~$0$. In other words, given feasibility, the optimal value of the MILP is always non-negative.
We next show that the optimal value is always positive if a feasible $\xb$ exists, for which the corresponding $\Ac(\xb)$ (containing $i$ by construction) leads to $\Kb^\ast \neq \zerob$ and a full-dimensional $\Rc^\ast$.
To this end, we recall that $\Ac(\xb)$ is determined by $\deltab^\ast$. We further note that the additional constraints~\eqref{eq:addConstraintsSaturation} imply $\xb^{(1)},\dots,\xb^{(m)} \in \Rc^\ast$ (analogously to $\xb \in \Rc^\ast$). Hence, we find
$$
\Ub_j^{(j)}-\Ub_j^\ast=\Kb_j^\ast \xb^{(j)}+\bb_j^\ast-\Kb_j^\ast \xb -\bb_j^\ast =\Kb_j^\ast (\xb^{(j)} - \xb)
$$ 
for every $j\in\{1,\dots,m\}$. Now, $\Kb^\ast \neq \zerob$ implies $\Kb_j^\ast \neq \zerob$ for at least one $j$. Due to $\Rc^\ast$ being full-dimensional, there exist $\xb^{(j)},\xb \in \Rc^\ast$ yielding a positive $\Ub_{j}^{(j)}-\Ub_{j}^\ast$. As a consequence, the optimal value in~\eqref{eq:costDeltaU} will be positive since the other terms in the cost function have already been shown to be non-negative (due to the $m$ independent $\xb^{(j)})$). Conversely, if the MILP returns $0$ as an optimal value, the active sets $\Ac(\xb)$ associated with feasible $\xb$ either correspond to $\Kb^\ast=\zerob$, lower dimensional $\Rc^\ast$, or both. Now, the former and the latter case are clearly irrelevant for computing the Lipschitz constant. However, also the  remaining case $\Kb^\ast\neq \zerob$ on some lower dimensional domain $\Rc^\ast$ is irrelevant if $\Fc_{\text{MPC}}$ is full-dimensional (as assumed). In fact, due to continuity of $\fb_{\text{MPC}}$ \cite[Thm.~4]{Bemporad2002}, the relevant gain $\Kb^\ast$ will then be captured by some neighboring segment on a full-dimensional domain.  
\end{proof}
\vspace{2mm}

Theorem~\ref{thm:saturatingInputs} provides another condition to potentially exclude the $i$-th constraint. Remarkably, the two conditions in the previous 
theorem and Lemma~\ref{lem:iConstraintInfeasible}, while conceptually different, are methodically closely related. To see this, note that the constraints of the corresponding MI problems only differ in terms of \eqref{eq:addConstraintsSaturation}. Now, it is easy to see that the constraints \eqref{eq:addConstraintsSaturation} are feasible whenever the corresponding constraints \eqref{eq:KKT_opt}--\eqref{eq:KKT_r} and \eqref{eq:KKT2MI} are feasible for $\xb$. Hence, we immediately find the following relation.

\begin{cor}
    The MI feasibilty problem in Lemma~\ref{lem:iConstraintInfeasible} is feasible if and only if the MILP in Theorem~\ref{thm:saturatingInputs} is feasible.
\end{cor}

As a consequence, one can only investigate the MILP in Theorem~\ref{thm:saturatingInputs} and exclude the $i$-th constraint if the MILP is either infeasible or returns $0$. 

\subsection{Exploiting symmetries in MPC}
\label{sec:symmetries}

Due to common symmetries in the constraints or the cost function, MPC often results in control laws, which likewise offer symmetries. Formalizing these symmetries can, e.g., be carried out analogously to \cite[Def.~1]{Danielson2015}. There, a symmetry is expressed in terms of invertible matrices $(\Thetab,\Omegab)$ satisfying
\begin{equation}
    \label{eq:symmetryDef}
    \Omegab \fb_{\text{MPC}}(\xb) = \fb_{\text{MPC}}(\Thetab \xb)
\end{equation}
for every $\xb \in \Fc_{\text{MPC}}$. Exploiting symmetries is, e.g., useful in the framework of explicit MPC \cite{Bemporad2002} since it allows to reduce the domain for which the explicit control law has to be computed (and stored). To specify this, we first note that multiple symmetries in terms of tuples $(\Thetab^{(1)},\Omegab^{(1)}),\dots,(\Thetab^{(\sigma)},\Omegab^{(\sigma)})$ can apply simultaneously with the canonical tuple $(\Ib_m,\Ib_n)$ being one of those. Then, we can substitute the constraint $\xb(0)\in \Xc$ in \eqref{eq:OCP} with $\xb(0)\in \Xc_{\text{fun}}$ for any choice of (the so-called fundamental domain) $\Xc_{\text{fun}}\subseteq \Xc$ satisfying
$$
 \Fc_{\text{MPC}} \subseteq \bigcup_{i=1}^\sigma \Thetab^{(i)} \Xc_{\text{fun}}.
$$
Clearly, in order to still enable the condensation to~\eqref{eq:QP}, it additionally makes sense to restrict ourselves to polyhedral sets $\Xc_{\text{fun}}$. Now, assuming for a moment that $\Xc$ and $\Xc_{\text{fun}}$ are characterized by the same number of hyperplanes. Then, it is easy to see that the substitution above does not alter the number of constraints in~\eqref{eq:QP}. Hence, it is not immediately clear how exploiting symmetries can be beneficial for our purposes. In this context, we first note that substituting $\Xc$ with a significantly smaller set $\Xc_{\text{fun}}$ (for the constraint associated with $\xb(0)$) often results in significantly more excluded constraints by the procedures related to Lemma~\ref{lem:iConstraintInfeasible} and Theorem~\ref{thm:saturatingInputs} (see our benchmark in Sect.~\ref{sec:Examples}). Moreover, symmetries often yield relations like $\Kb^{(i)}=-\Kb^{(j)}$ implying  ${\|\Kb^{(i)}\|_p=\|\Kb^{(j)}\|_p}$. However, symmetries do not always result in such trivial relations, in particular, in light of norms. To see this, note that~\eqref{eq:symmetryDef} in combination with the structure~\eqref{eq:fPWA} provides relations~like
$$
 \Kb^{(i)} \xb + \bb^{(i)}  = \Omegab^{-1} \Kb^{(j)} \Thetab \xb + \Omegab^{-1}\bb^{(j)}.
$$
Hence, instead of evaluating $\fb_{\text{MPC}}$ for some ${\xb \in \Rc^{(i)}}$, we could also make use of segment $j$ containing $\Thetab \xb$. This would allow us to skip segment $i$ in the context of the Lipschitz constant computation and to consider $\|  \Omegab^{-1} \Kb^{(j)} \Thetab \|_p$ instead. This observation can be exploited in two ways. 
First, we could simply evaluate $\sigma$ instances of the final MILP resulting for the tightened constraint $\xb(0)\in \Xc_{\text{fun}}$ in order to capture all transformed segments via $\|  \Omegab^{-1} \Kb^{(j)} \Thetab \|_p$ for every of the $\sigma$ tuples $(\Thetab^{(i)},\Omegab^{(i)})$. Second and more efficiently, we can evaluate the final MILP only once and consider only those transformations resulting in invariant norms, i.e.,
\begin{equation}
    \label{eq:invariantNormTransformation}
    \|\Omegab^{-1} \Kb \Thetab \|_p = \| \Kb \|_p \quad \text{for every}\,\,\, \Kb\in \R^{m\times n}.
\end{equation}
A sufficient condition for such transformations is as follows. 

\begin{lem}
\label{lem:normOmegaTheta}
    Let $\|\Omegab\|_p=1$ and $\|\Thetab\|_p=1$, then~\eqref{eq:invariantNormTransformation} holds.
\end{lem}

\begin{proof}
 We first note that invertability of $\Omegab$ and $\Thetab$ implies  $\|\Omegab^{-1}\|_p=\|\Omegab\|_p^{-1}=1$  and ${\|\Thetab^{-1}\|_p=1}$. Hence,
     $$
\|\Omegab^{-1} \Kb \Thetab \|_p \leq \|\Omegab^{-1}\|_p \| \Kb \|_p  \|\Thetab \|_p = \| \Kb \|_p
$$
due to sub-multiplicativity. On the other hand, 
\begin{align*}
   \| \Kb \|_p &= \| \Omegab \Omegab^{-1} \Kb \Thetab \Thetab^{-1} \|_p \\
   &\leq \| \Omegab \|_p  \| \Omegab^{-1} \Kb \Thetab \|_p \| \Thetab^{-1} \|_p =  \| \Omegab^{-1} \Kb \Thetab \|_p. 
\end{align*}
In combination, an inclusion results, which proves \eqref{eq:invariantNormTransformation}.
\end{proof}
\vspace{2mm}

While restrictive, common symmetries in MPC often satisfy the conditions in Lemma~\ref{lem:normOmegaTheta} (see, e.g., the examples in Sect.~\ref{sec:Examples}). It remains to comment on the identification of symmetries. 
In this context, we refer to the methods from \cite{Danielson2015}, which allow to identify tuples $(\Thetab,\Omegab)$ that reflect a symmetry purely based on $\Ab$, $\Bb$, $\Pb$, $\Qb$, $\Rb$, $\Xc$, and $\Uc$, i.e., without computing the explicit control law.

\subsection{Combined approaches}
\label{subsec:combinedApproaches}

We are ready to combine our approaches for a more efficient computation of Lipschitz constants. 
As already indicated, various combinations of the proposed tools can be considered. We specify two variants that will be used for the numerical benchmark in Section~\ref{sec:Examples}.
The variants differ in whether symmetries are exploited (according to the previous section) or not. Hence, slightly neglecting the additional effort for the identification of symmetries and a fundamental domain $\Xc_{\text{fun}}$, the crucial difference is that either $\xb(0) \in \Xc$ or $\xb(0) \in \Xc_{\text{fun}}$ is considered as a constraint for the initial state in~\eqref{eq:OCP}. Apart from this difference, all following steps are identical. In fact, we first condense the OCP to a QP of the form~\eqref{eq:QP}. We then use the MILP in Theorem~\ref{thm:saturatingInputs} to reduce the number of constraints. More precisely, we investigate for each constraint $i$ whether the MILP is infeasible or offers the optimal objective value $0$. In any of these cases, the $i$-th constraint is deleted. For simplicity of notation, we do not introduce different instances of the QP parameters for the two variants or during the constraint reduction. In fact, we simply assume that the previous instances are overwritten. Once the reduced QP is obtained, we solve the following MILP in order to compute the Lipschitz constant $L_1^\ast=l^\ast$:
\begin{equation}
    \label{eq:finalMILP}
 \max_{\substack{\xb, l,\,\Kb^\ast,  \Kb^+,\, \Kb^-,\,\cb,\,\Ub^\ast, \Ub^{(1)}, \dots, \Ub^{(n)},\\ \deltab^\ast, \deltab^{(1)},\dots,\deltab^{(m+1)}, \lambdab^\ast, \lambdab^{(1)},\dots,\lambdab^{(n)}, \rb^\ast, \rb^{(1)},\dots,\rb^{(n)}}} l 
\end{equation}
subject to \eqref{eq:KKT_opt}--\eqref{eq:KKT_r}, \eqref{eq:KKT2MI}, \eqref{eq:augmentedConstraints}, \eqref{eq:KfromUs}, and \eqref{eq:newNormConditions}.
Despite the MILP formulation of the reduced QP, a central element is the novel norm computation according to Lemma~\ref{lem:normMIFeas}.
Again, $L_\infty^\ast$ can be computed analogously.

\section{Numerical benchmark}
\label{sec:Examples}

\begin{table*}[h]
    \caption{Example systems from the literature.}
    \vspace{4mm}
    \centering
    \begin{tabular}{cccllccccl}
        \toprule
        No. & $\Ab$ & $\Bb$ & $\Xc$ & $\Uc$ & $\Qb$ & $\Rb$ & $N$ & $\Xc_{\text{fun}}$ & Reference \\  
        \midrule
        1. & $\begin{pmatrix}
            1 & 1 \\ 
            0 & 1
        \end{pmatrix}$ 
        & $\begin{pmatrix}
            0.5 \\
            1 
        \end{pmatrix}$ 
        & $\begin{array}{l}
             |\xb_1| \leq 25  \\
             |\xb_2| \leq 5
        \end{array}$ 
        & $\begin{array}{l}
             |\ub| \leq 1  
        \end{array}$ & $\Ib$ & $0.1$ & $10$ &
        $\begin{array}{c}
            |\xb_1| \leq 25 \!\!\!\!\!\!\!\!\!\!\!\!  \\
             0 \leq \xb_2 \leq 5
        \end{array}$ 
        & \cite[Eqs.~(2.8)--(2.9)]{Gutman1987} \\[4ex]

        2. & $\begin{pmatrix}
            0 & 1 \\ 
            1 & 0
        \end{pmatrix}$ 
        & $\begin{pmatrix}
            2 \\
            4 
        \end{pmatrix}$ 
        & $\begin{array}{l}
             |\xb_1| \leq 5  \\
             |\xb_2| \leq 5
        \end{array}$ 
        & $\begin{array}{l}
             |\ub| \leq 1  
        \end{array}$ & $\Ib$ & $4.5$ & $8$ 
        & 
        $\begin{array}{c}
            |\xb_1| \leq 5 \!\!\!\!\!\!\!\!\!  \\
             0 \leq \xb_2 \leq 5
        \end{array}$ 
        & \cite[Ex.~3]{SchulzeDarup2016_ECC_MPC} \\[4ex]
        
        3. & $\begin{pmatrix}
            1.1 & 0.2 \\ 
            -0.2 & 1.1
        \end{pmatrix}$ 
        & $\begin{pmatrix}
            0.5 & 0 \\
            0 & 0.4
        \end{pmatrix}$ 
        & $\begin{array}{l}
             |\xb_1| \leq 5  \\
             |\xb_2| \leq 5 
        \end{array}$ 
        & $\begin{array}{l}
             |\ub_1| \leq 1  \\
             |\ub_2| \leq 1 
        \end{array}$ & $\Ib$ & $0.1 \Ib$ & $3$ 
        & 
        $\begin{array}{c}
            |\xb_1| \leq 5 \!\!\!\!\!\!\!\!\!  \\
             0 \leq \xb_2 \leq 5
        \end{array}$ 
        & \cite[Ex.~2.26]{SchulzeDarup2014} \\[4ex]

        4. & $\begin{pmatrix}
            2 & 0 \\ 
            0 & 2
        \end{pmatrix}$ 
        & $\begin{pmatrix}
            1 & 1 \\
            1 & -1 
        \end{pmatrix}$ 
        & $\begin{array}{l}
             |\xb_1| \leq 5  \\
             |\xb_2| \leq 5
        \end{array}$ 
        & $\begin{array}{l}
             |\ub_1| \leq 1 \\
             |\ub_2| \leq 1
        \end{array}$ & $\Ib$ & $\Ib$ & $10$ 
        & 
        $\begin{array}{c}
             0 \leq \xb_1 \leq 5 \\
             0 \leq \xb_2 \leq 5
        \end{array}$ 
        & \cite[Ex.~1]{Danielson2015} \\[4ex]
        
        5. & $\begin{pmatrix}
            1 & 0.5 & 0.125 \\ 
            0 & 1 & 0.5 \\
            0 & 0 & 1
        \end{pmatrix}$ 
        & $\begin{pmatrix}
            0.02 \\
            0.125 \\
            0.5
        \end{pmatrix}$ 
        & $\begin{array}{l}
             |\xb_1| \leq 20  \\
             |\xb_2| \leq 3 \\
             |\xb_3| \leq 1
        \end{array}$ 
        & $\begin{array}{l}
             |\ub| \leq 0.5  
        \end{array}$ & $\Ib$ & $1$ & $3$ 
        & 
        $\begin{array}{c}
            |\xb_1| \leq 20 \!\!\!\!\!\!\!\!\!\!\!  \\
            |\xb_2| \leq 3 \!\!\!\!\!\!\!\! \\
            0 \leq \xb_3 \leq 1
        \end{array}$ 
        & \cite[Rem.~4.8]{Gutman1987} \\
        \bottomrule
        \label{tab:exampleSys}
    \end{tabular}
\end{table*}

\begin{table*}[h!]
    \caption{Computation of Lipschitz constants for different systems.}
    \vspace{4mm}
    \centering
    \begin{tabular}{ccrrrrrrrrrrrr}
        \toprule
        \multicolumn{3}{c}{} & \multicolumn{5}{c}{MILP from \cite[Thm.~5.3]{Fabiani2022}} & \multicolumn{5}{c}{MILP \eqref{eq:finalMILP} using Lemma \ref{lem:normMIFeas}, Theorem \ref{thm:saturatingInputs}, and symmetries} \\
        \cmidrule(lr{.75em}){4-8} \cmidrule(lr{.75em}){9-13} 
        No.\!\!\!\! & $\Tc\!\!\!\!$ & $\#_{\Rc}$\!\!\!\! & $\#_{\deltab}\!\!\!\!$ & $L^\ast_1$ & Time $[s]$\!\!\!\! & $L^\ast_{\infty}$ & Time $[s]$\!\!\!\! & $\#^{(1)}_{\deltab}\!\!\!\!$ & $\#^{(2)}_{\deltab}\!\!\!\!$ & Time for $L^*_1$ $[s]$\!\!\!\! & Time for $L^*_{\infty}$ $[s]$\!\!\!\! & Preprocessing $[s]$ \\  
        \midrule
        1.\!\!\!\! & $\Xc$\!\! & $211$ & $64$ & $1.89$ & $103.72$ & $1.27$ & $55.51$ & $22$ & $12$ & $0.1735$ & $0.2481$ & $0.0115$ \\[1ex]
        1.\!\!\!\! & $\Sc$\!\! & $195$ & $70$ & $1.89$ & $9.98$ & $1.27$ & $43.90$ & $22$ &  $12$ & $0.1453$ & $0.2289$ & $0.0156$ \\[2ex]
        2.\!\!\!\! & $\Xc$\!\! & $35$ & $52$ & $0.50$ & $15.24$ & $0.50$ & $32.44$ & $32$ & $17$ & $0.6331$ & $0.4759$ & $0.0143$ \\[1ex]
        2.\!\!\!\! & $\Sc$\!\! & $35$ & $58$ & $0.50$ & $19.05$ & $0.50$ & $43.24$ & $32$ & $17$ & $0.3275$ & $0.5374$ & $0.0136$ \\[2ex]
        3.\!\!\!\! & $\Xc$\!\! & $79$ & $28$ & $16.10$ & $2.59$ & $11.70$ & $3.37$ & $18$ & $15$ & $0.3079$ & $0.5262$ & $0.0587$ \\[1ex]
        3.\!\!\!\! & $\Sc$\!\! & $93$ & $32$ & $20.13$ & $1.75$ & $14.63$ & $1.40$ & $16$ & $16$ & $1.2126$ & $1.6429$ & $0.0863$ \\[2ex]
        4.\!\!\!\! & $\Xc$\!\! & $491$ & $84$ & $1.69$ & $71.36$ & $1.69$ & $242.74$ & $36$ & $12$ & $0.2590$ & $0.2520$ & $0.0145$ \\[1ex]
        4.\!\!\!\! & $\Sc$\!\! & $441$ & $88$ & $1.69$ & $66.14$ & $1.69$ & $85.64$ & $36$ & $12$ & $1.6266$ & $0.1910$ & $0.0248$ \\[2ex]
        5.\!\!\!\! & $\Xc$\!\! & $117$ & $30$ & $12.00$ & $2.11$ & $8.00$ & $3.49$ & $28$ & $22$ & $0.3966$ & $0.9440$ & $0.0434$ \\[1ex]
        5.\!\!\!\! & $\Sc$\!\! & $107$ & $46$ & $12.00$ & $9.69$ & $8.00$ & $22.87$ & $28$ & $22$ & $1.1988$ & $2.4368$ & $0.0539$ \\
        \bottomrule
        \label{tab:results}
    \end{tabular}
\end{table*}

We demonstrate the effectiveness of the proposed procedures by applying them to examples from the literature summarized in Table~\ref{tab:exampleSys}. All MILP in this section are solved using the mixed-integer solver from \cite{mosek} with constants set to $\overline{\lambdab}=\overline{\rb}=M=10^4$. For every example, we first compute the minimum Lipschitz constants $L^\ast_1$ and $L^\ast_\infty$ according to the method in \cite[Thm.~5.3]{Fabiani2022}, which also served as the starting point for our investigations. The required computation time is listed in the sixth and eighth column of Table~\ref{tab:results}, respectively.

Next, we follow the two variants in Section \ref{subsec:combinedApproaches} in order to apply our novel procedures. Regarding the variant exploiting symmetries, we note that all examples in Table~\ref{tab:exampleSys} offer rotational symmetries. More precisely, the matrices $\Thetab$ and $\Omegab$ are of the form
\begin{equation*}
    \Thetab=\begin{pmatrix}
        \cos(\varphi) & \sin(\varphi) \\
        -\sin(\varphi) & \cos(\varphi)
    \end{pmatrix}, 
\end{equation*}
with $\varphi\in\{\pi,2\pi\}$ and $\Omegab\in\{-1,1\}$ for Systems~$1$ and~$2$. For the two following systems, 
we have $\Thetab$ as before and 
\begin{equation*}
\Omegab=\begin{pmatrix}
        \cos(\varphi) & -\sin(\varphi) \\
        \sin(\varphi) & \cos(\varphi)
    \end{pmatrix}=\Thetab^\top
\end{equation*}
with $\varphi\in\{\pi,2\pi\}$ for System $3$ and $\varphi\in\{\nicefrac{\pi}{2},\pi,\nicefrac{3\pi}{2},2\pi\}$ for System $4$. For System $5$, $\Thetab$ likewise reflects rotation matrices (which we omit for brevity) and $\Omegab\in\{-1,1\}$. 
In all cases, we easily verify that the condition in Lemma~\ref{lem:normOmegaTheta} (i.e., $||\Omegab||_p = ||\Thetab||_p=1$) holds. Hence, we can restrict our analysis to the fundamental domains in the ninth column of Table~\ref{tab:exampleSys} without the need to consider multiple instances of the MILP~\eqref{eq:finalMILP}. Before analyzing the performance of our procedures, we note that we consider two different choices for the terminal set $\Tc$ for each example. First, we simply choose $\Tc=\Xc$ (i.e., no terminal constraints). Second, we consider $\Tc=\Sc$ with $\Sc$ denoting the largest positively invariant set, where the linear quadratic regulator (LQR) can be applied without violating constraints. The computation has been carried out according to \cite{Gilbert1991}.

Now, in Table~\ref{tab:results}, we list the number of binary variables of the MILP in \cite[Thm.~5.3]{Fabiani2022} in the column $\#_{\deltab}$. Further,  we also list the number of regions $\Rc^{(i)}$ of the explicit control law, computed using the multi-parametric toolbox (MPT, \cite{MPT3}), under $\#_{\Rc}$ as an orientation. Key performance indicators of our procedures are listed in columns nine to 13. First, $\#^{(1)}_{\hat{\deltab}}$ is the number of binary variables for the simplified MILP without considering symmetries. Second, $\#^{(2)}_{\hat{\deltab}}$ reflects the same figure with symmetries. Next, the total computation times for evaluating $L_1^\ast$ and $L_\infty^\ast$ are listed, respectively. These include the times for all preprocessing steps (such as constraint elimination), which are exclusively listed in the last column for completeness.  

As apparent from Table~\ref{tab:results}, the time required to compute the Lipschitz constants can be significantly reduced for all examples. Moreover, one can observe that the computation time for the preprocessing is negligible compared to the time required to compute $L^\ast_p$. In fact, although MILP are solved during the preprocessing, the overall time can be reduced in all cases. Furthermore, for the proposed method, the variation of the computation times between the examples is lower. This may indicate, that the proposed method scales better with the model complexity.    

\begin{figure}[ht]
\begin{minipage}[h]{0.48\columnwidth}
\centering
\includegraphics[trim={0cm 0cm 0cm 0cm},clip,scale=0.5]{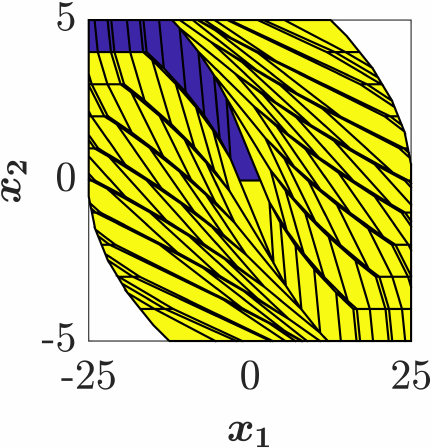}
\caption{Sate space partition\\ of the control law.}
\label{fig:OCLPartition}
\end{minipage}
\begin{minipage}[h]{0.48\columnwidth}
\centering
\includegraphics[trim={0cm 0cm 0cm 0cm},clip,scale=0.5]{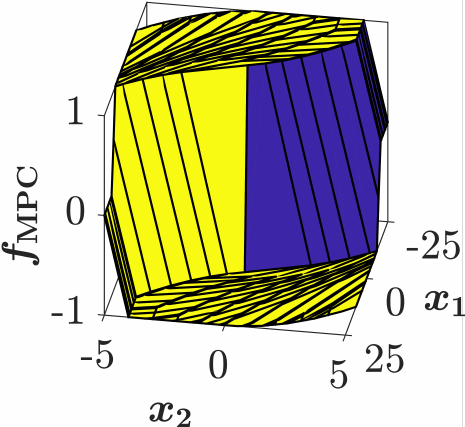}
\caption{Optimal control law for System $1$ with $\Tc=\Xc$.}
\label{fig:OCL}
\end{minipage}
\end{figure}

Finally, Figure~\ref{fig:OCLPartition} highlights in purple the regions that are considered by the MILP~\eqref{eq:finalMILP} for computing the Lipschitz constant after reducing the number of binary variables according to Lemma~\ref{lem:iConstraintInfeasible} and Theorem~\ref{thm:saturatingInputs} under consideration of symmetries. As apparent from Figure~\ref{fig:OCL}, Theorem~\ref{thm:saturatingInputs} allows to identify all regions with a local gain of zero. Moreover, all regions that have the same local gain due to symmetries are also excluded from the computation. Thus only regions that are relevant for computing a Lipschitz constant are considered.
\section{Conclusion}\label{sec:Conclusion}
We presented an efficient set of methods to compute the minimum Lipschitz constant of an MPC law. The method adopts a known MILP for the computation of Lipschitz constants and uses various procedures (see Lemma~\ref{lem:normMIFeas}, Lemma~\ref{lem:iConstraintInfeasible} and Theorem~\ref{thm:saturatingInputs}) to reduce the number of binary variables in the MILP and thus its complexity. The most powerful reduction builds on exploiting saturation and symmetries of the control law. However, the proposed reduction steps can also be applied to systems without symmetries. 
This allows an efficient computation of the minimum Lipschitz constant even for complex systems of moderate size. For future research it would be interesting to investigate if the MILP~\eqref{eq:finalMILP} can be adopted for the computation of Lipschitz constants of the piecewise quadratic optimal value function $V(\xb)$ of the OP~\eqref{eq:OCP}.


\end{document}